\documentclass[11pt,a4paper]{amsart}
\usepackage{mathrsfs}
\usepackage{syntonly}
\usepackage{amsmath}
\usepackage{amsthm}
\usepackage{amsfonts}
\usepackage{amssymb}
\usepackage{latexsym}
\usepackage{amscd,amssymb,amsopn,amsmath,amsthm,graphics,amsfonts,mathrsfs,accents,enumerate,verbatim,calc}
\usepackage[dvips]{graphicx}
\usepackage[colorlinks=true,linkcolor=red,citecolor=blue]{hyperref}
\usepackage[all]{xy}
\usepackage{enumitem}

\date{}
\allowdisplaybreaks[4] \footskip=15pt
\renewcommand{\uppercasenonmath}[1]{}

\numberwithin{equation}{section} \theoremstyle{plain}
\newtheorem{lem}{Lemma}[section]
\newtheorem{cor}[lem]{Corollary}
\newtheorem{prop}[lem]{Proposition}
\newtheorem{thm}[lem]{Theorem}

\newtheorem{Ex}[lem]{Example}
\newtheorem{Quest}[lem]{Question}
\newtheorem{Property}[lem]{Property}
\newtheorem{Properties}[lem]{Properties}
\newtheorem{Subprops}{}[lem]
\newtheorem{Para}[lem]{}

\newtheorem{rem}[lem]{Remark}

\newtheorem*{ack*}{ACKNOWLEDGEMENTS}

\newcommand{\pf}{\noindent\begin {proof}}
\newcommand{\epf}{\end{proof}}

\newcommand{\X}{\mathcal{X}}

\newcommand{\Y}{\mathcal{Y}}

\newcommand{\A}{\mathcal{A}}

\newcommand{\Arr}{{\rm Arr}}

\newcommand{\ME}{{\rm ME}}

\newcommand{\Ab}{{\rm {\bf \mbox{Ab}}}}

\newcommand{\Ext}{{\rm Ext}}

\newcommand{\Hom}{{\rm Hom}}

\newcommand{\op}{{\rm {\footnotesize op}}}

\newcommand{\mcA}{\mathcal{A}}

\newcommand{\mcE}{\mathcal{E}}

\newcommand{\mcI}{\mathcal{I}}
\newcommand{\mcJ}{\mathcal{J}}

\newcommand{\mcM}{\mathcal{M}}

\newcommand{\mcX}{\mathcal{X}}
\newcommand{\mcY}{\mathcal{Y}}

\begin{document}
\begin{center}
{ \bf \MakeUppercase{Cotorsion pairs and Enochs Conjecture for object ideals}}

\vspace{0.5cm}  Dandan Sun, Qikai Wang and Haiyan Zhu\footnote{Corresponding author.\\
\indent Haiyan Zhu was supported by the NSF of China (12271481).}
\end{center}
\bigskip
\centerline { \bf  Abstract}
\leftskip10truemm \rightskip10truemm \noindent
Let $\mcI$ and $\mcJ$ be object ideals in an exact category $(\mcA; \mcE)$. It is proved that $(\mcI,\mcJ)$ is a perfect ideal cotorsion pair if and only if $({\rm Ob}(\mcI),{\rm Ob}(\mcJ))$ is a perfect cotorsion pair, where ${\rm Ob}(\mcI)$ and ${\rm Ob}(\mcJ)$ is the objects of $\mcI$ and $\mcJ$, respectively. If in addition $(\mcA; \mcE)$ has enough projective objects and injective objects, and $\mcJ$ is enveloping, then  $(\mcI,\mcJ)$ is a complete ideal cotorsion pair if and only if $({\rm Ob}(\mcI),{\rm Ob}(\mcJ))$ is a complete cotorsion pair. This gives a partial answer to the question posed by Fu, Guil Asensio, Herzog and Torrecillas. Moreover, for any object ideal $\mcI$ in the category of left $R$-modules, it is proved that $\mcI$ satisfies Enochs Conjecture if and only if ${\rm Ob}(\mcI)$ satisfies Enochs Conjecture. Applications are given to projective morphisms and  ideal cotorsion pairs $(\mcI,\mcJ)$ of object ideals under certain conditions.
\leftskip10truemm \rightskip10truemm \noindent
\\[2mm]
{\bf Keywords:} Object ideals; Enochs Conjecture; complete ideal cotorsion pairs;  perfect ideal cotorsion pairs.\\
{\bf 2020 Mathematics Subject Classification:} 18G10, 18G25, 16D90.

\leftskip0truemm \rightskip0truemm

\section{Introduction}
The notion of a cotorsion pair was originally defined in the category of abelian groups \cite{SL}. It was later extended to abelian categories and exact categories. In recent years, the study of cotorsion pairs has gained significance in the investigation of precovers and preenvelopes (approximation theory), particularly in the context of proving the flat cover conjecture \cite{BEE}. In short, a cotorsion pair in an exact category $(\mcA; \mcE)$ is a pair $(\X,\Y)$ of classes of objects of $\A$, such that each of which is the orthogonal complement of the other with respect to the Yoneda  bifunctor functor $\Ext(-,-)$. If in addition every object in $\A$ has a special $\X$-cover and a special $\Y$-envelope, then it is  said to be a \emph{complete} cotorsion pair.
Following \cite[Definition 2.3.1]{Gobel}, the cotorsion pair $(\X,\Y)$ is called \emph{perfect} if every object in $\A$ has an $\X$-cover and a $\Y$-envelope. The term ``perfect" comes from the classical result of Bass \cite{Bass1960} characterizing left perfect ring by the property of perfectness of the cotorsion ($\textrm{Proj}$$R$, $R$-Mod), where $R$-Mod is the class of left $R$-modules and $\textrm{Proj}$$R$ denotes the classes of projective left $R$-modules.

In an abstract category, objects and morphisms are two essential components; and by a well-known embedding from a category to its morphism category, objects can be viewed as
special morphisms. In order to give morphisms and ideals of categories equal importance, Fu, Guil Asensio, Herzog and Torrecillas \cite{FGHT} introduced ideal cotorsion pairs and the ideal approximation theory. In this theory, the role of the objects and subcategories in classical approximation theory is replaced by morphisms and ideals of the category.
If $\mcI$ is an ideal of $\mcA$, define ${\mcI}^{\perp}$ to be the ideal of morphisms $j$ such that $\Ext(i,j)=0$ for every $i\in{\mcI}$, and define $^{\perp}{\mcI}$ dually. An \emph{ideal cotorsion pair} in $(\mcA; \mcE)$ is a pair $(\mcI,\mcJ)$ of ideals of $\mcA$ satisfying $\mcJ={{\mcI}^{\perp}}$ and $\mcI={^{\perp}{\mcJ}}$. Ideal approximation theory began with an ideal version of Salce's Lemma [7, Theorem 18], which states that if $(\mcI,\mcJ)$ is an ideal cotorsion pair in an exact category $(\mcA; \mcE)$ with enough projective
morphisms and injective morphisms, then every object in $\mcA$ admits a special $\mcJ$-preenvelope if and only if every object has a special $\mcI$-precover. Such cotorsion pairs $(\mcI,\mcJ)$ are called \emph{complete}. Following \cite{Gobel}, the ideal cotorsion pair $(\mcI,\mcJ)$ is called \emph{perfect} if every object in $\mcA$ admits a $\mcJ$-envelope and an $\mcI$-cover.

The classical situation of a complete cotorsion pair $(\mcX, \mcY)$ in an exact category $(\mcA; \mcE)$
yields an example of a complete ideal cotorsion pair, given by the pair $(\mcI(\mcX), \mcI(\mcY))$
of the corresponding object ideals, where $\mcI(\mcX)$ (resp., $\mcI(\mcY)$) denotes the ideal of those morphisms that factor through an object of the
additive subcategory $\mcX\subseteq \mcA$ (resp., $\mcY\subseteq \mcA$). However, the converse problem, due to Fu, Asensio, Herzog and Torrecillas \cite[Question 29]{FH},
that is if the cotorsion pair $({\rm Ob}(\mcI),{\rm Ob}(\mcJ))$ is necessarily complete provided that $\mcI$ and $\mcJ$ are object ideals in $\mcA$ such that $(\mcI,\mcJ)$ is a complete ideal cotorsion pair, is still open. Here ${\rm Ob}(\mcI):=\{X\in{\mcA}~|~1_{X}\in{\mcI}\}$ (resp., ${\rm Ob}(\mcJ):=\{X\in{\mcA}~|~1_{X}\in{\mcJ}\}$) is a subcategory of $\mcA$. The first aim of this note is to address this problem under certain conditions.

An interesting and deep theory of Enochs says that any precovering class of modules closed under direct limits is covering (see \cite{Enochs1981}). The converse problem,
that is if a covering class of modules is necessarily closed under direct limits, is still open and
is known as \emph{Enochs Conjecture} for objects (see \cite[Open problems 5.4]{Gobel}).
Some significant advancements have been made towards the solution of this conjecture in recent years (see \cite{AST-2018,BL,Bazzoni,BPS,Bazzoni-Saroch,Saroch,Saroch2} for instance). 
On the other hand, Estrada, Guil Asensio and Odabasi proved that any precovering class of morphisms in the category $R$-Mod of left $R$-modules closed under direct limits is covering (see \cite[Proposition 3.9]{EGO}). It seems natural to ask whether each covering ideal in $R$-Mod is closed under direct limits. A hypothetical general positive answer to this question is called ``\emph{Enochs Conjecture}" for morphisms.
The second aim of this note is to establish the relations between Enochs Conjecture for objects and morphisms.

We now outline the results of the paper. In Section 2, we summarize some preliminaries and
basic facts which will be used throughout the paper.

In Section 3, for any object ideal $\mcI$ in an exact category $(\mcA; \mcE)$, we first show that $\mcI$ is covering (resp., enveloping) if and only if ${\rm Ob}(\mcI)$ is covering (resp., enveloping) (see Proposition \ref{thm:covering-object}). Then we will prove that $(\mcI,\mcJ)$ is a perfect ideal cotorsion pair if and only if $({\rm Ob}(\mcI),{\rm Ob}(\mcJ))$ is a perfect cotorsion pair whenever $\mcI$ and $\mcJ$ are object ideals in $\mcA$ (see Theorem \ref{cor3}(1)). If in addition $(\mcA; \mcE)$ has enough projective objects and injective objects, and $\mcJ$ is enveloping, then  $(\mcI,\mcJ)$ is a complete ideal cotorsion pair if and only if $({\rm Ob}(\mcI),{\rm Ob}(\mcJ))$ is a complete cotorsion pair (see Theorem \ref{cor3}(2)). This gives a partial answer to Question 29 in \cite{FGHT}.

In Section 4, for any object ideal $\mcI$ in $R$-Mod, it is proved that $\mcI$ satisfies Enochs Conjecture if and only if ${\rm Ob}(\mcI)$ satisfies Enochs Conjecture (see Theorem \ref{cor1}). As a consequence, it is shown that the class of projective morphisms in $R$-Mod satisfies Enochs Conjecture, which can be seen as the ideal version of the result of Bass \cite{Bass1960} that the class of projective objects in $R$-Mod satisfies Enochs Conjecture. As another consequence, for any
object ideals $\mcI$ and $\mcJ$ in ${R}$-{\rm Mod}, it is proved in Corollary \ref{cor:4.4} that $\mcI$ satisfies Enochs Conjecture provided that $(\mcI,\mcJ)$ is an ideal cotorsion pair such that $\mcJ$ is closed under direct limits. Some conditions are given to show that the validity of Enchos Conjecture for objects of a general ideal $\mcI$ (not necessarily object ideal) implies the validity of Enchos Conjecture of ${\rm Ob}(\mcI)$ (see Propositions \ref{prop4.5} and \ref{prop-last}).

\section{Preliminary}

The assumptions, the notation, and the definitions from this section will be used throughout the
paper.  For more details the reader can consult \cite{EJ2,FGHT,FH,Gobel}.

\subsection{Exact categories}
Recall from \cite{TBh10,Keller,Quillen} that an \emph{exact category} is a pair $(\mcA; \mcE),$  where $\mcA$ is an additive category and $\mcE$ is a class of ``short exact sequences": That is,
an actual kernel-cokernel pair 
$$\Xi \colon \xymatrix@1@C=45pt{A \ar[r] & B \ar[r] & C}.$$ In what follows, we call such a sequence a \emph{conflation}, and call $A\to B$ (resp., $B\to C$) an \emph{inflation} (resp., \emph{deflation}). Many authors use the alternate terms \emph{admissible exact sequence}, \emph{admissible monomorphism} and \emph{admissible epimorphism}.
The class $\mathcal{E}$ of conflations must satisfy exact axioms, for details, we refer the reader to \cite[Definition 2.1]{TBh10}, which are inspired by the properties of short exact sequences in any abelian category. We often write $\mcA$ instead of $(\mcA; \mcE)$ when we consider only one exact structure on $\mathcal{A}$. 

%

Recall that an object $P$ in $\mcA$ is called \emph{projective} provided that any admissible epimorphism ending at $P$ splits.
The exact category $\mcA$ is said to have \emph{enough projective objects} provided that each object $X$ fits into a deflation $d: P \rightarrow X$ with $P$ projective. Dually one has the notions of injective
objects and enough injective objects.
%
%
%

\subsection{Ideals}
An {\em ideal}\;\;$\mcI$ of an additive category $\mcA$ is an additive subbifunctor of $\Hom : \mcA^{\op} \times \mcA \to \Ab;$ it associates to a pair $(A,B)$ of morphisms in $\mcA$ a subgroup 
$\mcI(A,B) \subseteq \Hom (A,B)$
and satisfies the usual conditions in the definition of an ideal in ring theory: for any morphism $g: A \to B$ in $\mcI,$ the composition $fgh$ belongs to $\mcI (X,Y),$ for any morphisms $f: B \to Y$ and $h: X \to A$ in $\mcA.$ 

If $A\in{\mcA}$ is an object, we say that $A$ belongs to an ideal $\mcI$ if the identity morphism $1_{A}$ belongs to
$\mcI(A,A),$ and we denote by $\textrm{Ob}(\mcI)\subseteq \mcA$ the full subcategory of objects of $\mcI$. In the other direction, every
additive subcategory $\mcX\subseteq{\mcA}$ gives rise to the ideal $\mcI(\mcX)$ generated by morphisms of the form $1_{X}$ for all $X\in{\mcX}$; this is the ideal of morphisms that factor through an object in $\mcX$. An ideal $\mcI$ of $\mcA$ is called an
\emph{object ideal} if it is generated by its objects, $\mcI = \mcI(\textrm{Ob}(\mcI))$ (see \cite{FGHT}).

\subsection{(Pre)covering and (pre)enveloping ideals}

Let $\mcI$ be a class of morphisims in an additive category $\mcA$. An $\mcI$-precover of an object $X\in{\mcA}$ is a morphism $f:A\rightarrow X$ in $\mcI$ such that any other morphism $f':A'\rightarrow X$ in $\mcI$ factors through $f$. It is said to be an \emph{$\mcI$-cover} if every $g$ with $fg=g$ is an isomorphism. 

In addition, if $\mcA$ is endowed with an exact structure $\mcE$, then $\mcI$-precover $f:A\rightarrow X$ of $X$ is said to be \emph{special} if $f$ is a deflation in $\mcE$ and it fits in a pushout diagram 
$$
\xymatrix@C=30pt@R=30pt{
K \ar[r] \ar[d]^{g} & A' \ar[r] \ar[d] & X \ar@{=}[d] \\
E \ar[r] & A \ar[r]^{f} & X}
$$
with $\Ext(i,g)=0$ for every $i\in{\mcI}$. The class $\mcI$ is said to be \emph{(special) (pre)covering} if every object in $\mcA$ has a (special) $\mcI$-(pre)cover. The notion of (special) $\mcI$-(pre)envelope is defined dually.

Let $\mcI$ be an ideal in $R$-Mod. Following \cite{AR,EG2014}, we say that $\mcI$ is \emph{closed under direct limits} if for any morphism $\{f_{i}:M_{i}\to N_{i}\}_{i\in{\mathbb{N}}}$ between directed systems of morphisms $\{g_{ij}:M_{i}\to M_{j}\}_{i\leqslant j}$ and $\{h_{ij}:N_{i}\to M_{j}\}_{i\leqslant j}$,
satisfying that $f_{i}\in{\mcI}$ for every $i\in{\mathbb{N}}$, the induced morphism $\lim\limits_{\longrightarrow}f:\lim\limits_{\longrightarrow}M_{i}\rightarrow \lim\limits_{\longrightarrow}N_{i}$ also belongs to $\mcI$. The following result is a direct consequence of \cite[Proposition 3.9]{EGO}.

\begin{prop}\label{prop:2.1} Let $R$ be a ring. Then any precovering class of morphisms in the category of left $R$-modules closed under direct limits is covering.
\end{prop}

\subsection{Ideal cotorsion pairs} As defined in the Introduction, two morphisms $a$ and $b$ in $\mcA$ are \emph{$\Ext$-orthogonal} if the morphism $\Ext(a,b) : \Ext(A_1, B_0) \to \Ext(A_0, B_1)$ of abelian groups is $0.$
If $\mcM$ is a class of morphisms in $\mcA,$ then ${^{\perp}}\mcM$ is the ideal of morphisms left $\Ext$-orthogonal to every $m \in \mcM.$ One defines a {\em right closed} ideal dually, and an \emph{ideal cotorsion pair} consists of ideals $(\mcI, \mcJ)$ for which $\mcI = {^{\perp}}\mcJ$ and $\mcJ^{\perp} = \mcI.$ Moreover,  an ideal cotorsion pair $(\mcI,\mcJ)$ is \emph{complete} if every object in $\mcA$ has a special $\mcJ$-preenvelope and a special $\mcI$-precover (see \cite{FGHT}). In the following, we call an ideal cotorsion pair $(\mcI,\mcJ)$ is \emph{perfect} if every object in $\mcA$ has a $\mcJ$-envelope and an $\mcI$-cover.

\subsection{ME-conflations}
The {\em arrow category} $\Arr (\mcA)$ of a category $\mcA$ is the category whose objects \linebreak $a: A_0 \to A_1$ are the morphisms (arrows) of $\mcA,$
and a morphism $f: a \to b$ in $\Arr (\mcA)$ is given by a pair of morphisms $f = (f_0, f_1)$ in $\mcA$ for which the diagram
$$
\xymatrix@C=30pt@R=30pt{A_0 \ar[r]^{f_0} \ar[d]^{a} & B_0 \ar[d]^b \\
A_1 \ar[r]^{f_1} & B_1}$$
commutes. In addition, the additive category $\Arr (\mcA)$ may be equipped with an exact structure for which the collection $\Arr (\mcE)$ of distinguished kernel-cokernel pairs
$\xi: \xymatrix@1{b \ar[r]^m & c \ar[r]^p & a,}$ is given by
morphisms of conflations in $(\mcA; \mcE),$
$$
\xymatrix@C=30pt@R=35pt{\Xi_0: B_0 \ar@<-2.5ex>[d]_{\xi} \ar@<2ex>[d]^b \ar[r]^-{m_0} & C_0 \ar[r]^{p_0} \ar[d]^c & A_0 \ar[d]^a \\
\Xi_1: B_1 \ar[r]^-{m_1} & C_1 \ar[r]^{p_1} & \; A_1.}$$

A conflation of arrows $\xi: \xymatrix@1@C=15pt{b \ar[r] & c \ar[r] & a}$ in the exact structure $(\Arr (\mcA); \Arr (\mcE))$ is said to be $\ME$ ({\em mono-epi}) if there exists a factorization
$$\xymatrix{
\Xi_0 : B_0 \ar[r] \ar@{=}@<2ex>[d] \ar@<-3ex>[dd]_{\xi} & C_0 \ar[r] \ar[d]^{c_1} & A_0 \ar[d]^a \\
 \;\; \;\;\;\; B_0 \ar[r] \ar@<2ex>[d]^b & C \ar[r] \ar[d]^{c_2} & A_1 \ar@{=}[d]  \\
\Xi_1 : B_1 \ar[r] & C_1 \ar[r] & A_1,
}$$ with $c = c_2 c_1.$ The notation $c = b \star a$ is used to indicate that the arrow $c$ is an $\ME$-extension of $a$ by $b.$

\section{\bf The completeness and perfectness for ideal cotorsion pairs of objects ideals}\label{application:Enochs-conjecture}

We begin with the following easy observation.
\begin{lem}\label{dir-sum-of-objec} Let $\mcI$ be an ideal in an exact category $(\mcA; \mcE)$. Then ${\rm Ob}(\mcI)$ is closed under direct summands.
\end{lem}
\begin{proof} Assume that $A\oplus B$ is in ${\rm Ob}(\mcI)$. Then ${\left(\begin{smallmatrix}1&0\\0&1\\\end{smallmatrix}\right)}:A\oplus B\to {A\oplus B}$ belongs to $\mcI$. Notice that $1_{A}={\left(\begin{smallmatrix}1\\0\\\end{smallmatrix}\right)}{\left(\begin{smallmatrix}1&0\\0&1\\\end{smallmatrix}\right)}
(1,0)$. Thus $1_{A}$ belongs to $\mcI$, as desired.
\end{proof}

The following lemma is a consequence of Remark 3.3 in \cite{EGO}.
\begin{lem}\label{prop:precovering-object}  Let $\mcI$ be an object ideal in an exact category $(\mcA; \mcE)$. 
 \begin{enumerate}
\item $\mcI$ is precovering (resp., preenveloping) if and only if ${\rm Ob}(\mcI)$ is precovering (resp., preenveloping);
\item Every ${\rm Ob}(\mcI)$-precover $f:X\to A$ of $A$ is also an $\mcI$-precover of $A$;
\item Every ${\rm Ob}(\mcI)$-preenvelope $f:A\to X$ of $A$ is also an $\mcI$-preenvelope of $A$. 
\end{enumerate}   
\end{lem}
\begin{proof} Thanks to \cite[Remark 3.3]{EGO}, it suffices to prove (2) and (3). Let  $f:X\to A$ be an ${\rm Ob}(\mcI)$-precover. If we choose any morphism $g:F\to A$ in the object ideal $\mcI$, then there exists an object $G$ in ${\rm Ob}(\mcI)$ with $\alpha:G\to A$ and $\beta:F\to G$ such that $g=\alpha\beta$. Since $f:X\to A$ is an ${\rm Ob}(\mcI)$-precover of $A$, there exists a morphism $\gamma:G\to X$ such that $\alpha=f\gamma$. This implies that $f(\gamma\beta)=(f\gamma)\beta=\alpha\beta=g$, and therefore $f:X\to A$ is an $\mcI$-precover. This finishes the
proof of (2), and the proof of (3) is similar.
\end{proof}

\begin{prop}\label{thm:covering-object} The following are true for any object ideal $\mcI$ in an exact category $(\mcA; \mcE)$.
 \begin{enumerate}
\item $\mcI$ is covering if and only if ${\rm Ob}(\mcI)$ is covering.

\item $\mcI$ is enveloping if and only if ${\rm Ob}(\mcI)$ is enveloping.

\end{enumerate}
\end{prop}
\begin{proof} We only prove (1), and (2) is dual to (1). For the ``only if" statement, we choose an $\mcI$-cover $f:X\to A$ for any object $A$ in $\mcA$. Next we claim that $X$ belongs to ${\rm Ob}(\mcI)$. By Lemma \ref{prop:precovering-object}(1), $A$ has an ${\rm Ob}(\mcI)$-precover $g:F\to A$ with $F$ in ${\rm Ob}(\mcI)$. Thus there exists $\beta:F\to X$ such that $g=f\beta$. On the other hand, it follows from Lemma \ref{prop:precovering-object}(2)  that there exists $\alpha:X\to F$ such that $f=g\alpha$. This yields that $f(\beta\alpha)=(f\beta)\alpha=g\alpha=f$. Because $f:X\to A$ is an $\mcI$-cover, we obtain that
$\beta\alpha$ is an isomorphism, and therefore $X$ is isomorphic to a direct summand of $F$. So $X$ belongs to ${\rm Ob}(\mcI)$ by Lemma \ref{dir-sum-of-objec}, as desired.

For the ``if" statement, there exists an ${\rm Ob}(\mcI)$-cover $f:X\to A$ for any object $A$ by hypothesis. Hence $f:X\to A$ is also an $\mcI$-cover, and the proof is similar to that of (2) in Lemma \ref{prop:precovering-object}.
\end{proof}

\begin{lem}\label{lem:easy} The following are true for any ideal cotorsion pair $(\mcI, \mcJ)$ in an exact category $(\mcA; \mcE)$:
 \begin{enumerate}
\item An object $X\in{{\rm Ob}(\mcI)}$ if and only if $\Ext_{\mcA} (X,j)=0$ for any $j\in{\mcJ}$.
\item An object $Y\in{{\rm Ob}(\mcJ)}$ if and only if $\Ext_{\mcA} (i,Y)=0$ for any $i\in{\mcI}$.
\end{enumerate}
\end{lem}
\begin{proof} The proof is straightforward by definition.
\end{proof}

\begin{lem}\label{lemma:3.3-1} Let $\mcI$ and $\mcJ$ be object ideals in an exact category $(\mcA; \mcE)$. Then $(\mcI,\mcJ)$ is an ideal cotorsion pair if and only if $({\rm Ob}(\mcI),{\rm Ob}(\mcJ))$ is a cotorsion pair.
\end{lem}
\begin{proof} Thanks to \cite[Theorem 2.8]{FGHT}, we only prove the ``only if" statement. Next we prove ${\rm Ob}(\mcJ)={{\rm Ob}(\mcI)^{\perp}}$, and the proof of ${\rm Ob}(\mcI)={^{\perp}{\rm Ob}(\mcJ)}$ is similar. By Lemma \ref{lem:easy}, it is straightforward to show the containment ${\rm Ob}(\mcJ)\subseteq{{\rm Ob}(\mcI)^{\perp}}$. For the reverse containment, we choose an object $M\in{{\rm Ob}(\mcI)^{\perp}}$. Let $f:X\rightarrow Y$ be an object ideal in $\mcI$. Then there exist an object $K\in{{\rm Ob}(\mcI)}$ and morphisms $g:K\to Y$ and $f:X\to K$ such that $f=gh$. Since $\Ext_{\mcA}(K,M)=0$ by hypothesis, and therefore $\Ext_{\mcA}(f,M)=\Ext_{\mcA}(h,M)\Ext_{\mcA}(g,M)=0$. Thus Lemma \ref{dir-sum-of-objec} yields that $M\in{{\rm Ob}(\mcJ)}$, as desired.
\end{proof}

We are now in a position to state and prove the main result of this section. It should be noted that Theorem \ref{cor3}(2) gives a partial answer to \cite[Question 29]{FGHT}.
\begin{thm}\label{cor3} The following are true for object ideals $\mcI$ and $\mcJ$ in an exact category $(\mcA; \mcE)$:
\begin{enumerate}
\item $(\mcI,\mcJ)$ is a perfect ideal cotorsion pair if and only if $({\rm Ob}(\mcI),{\rm Ob}(\mcJ))$ is a prefect cotorsion pair.
\item If $(\mcA; \mcE)$ has enough projective objects and injective objects, and $\mcJ$ is enveloping, then $(\mcI,\mcJ)$ is a complete ideal cotorsion pair if and only if $({\rm Ob}(\mcI),{\rm Ob}(\mcJ))$ is a complete cotorsion pair.
\end{enumerate}
\end{thm}
\begin{proof} We only prove (2), and (1) is a direct consequence of Proposition \ref{thm:covering-object} and Lemma \ref{lemma:3.3-1}. By \cite[Theorem 2.8]{FGHT}, it suffices to show the ``only if" statement. Let $A$ be an object in $\mcA$. Then there exists a conflation
$$\xymatrix@1@C=45pt{K \ar[r] & P \ar[r] & A,}$$
where $P$ is a projective object in $\mcA$. Since $\mcJ$ is enveloping, it follows from Proposition \ref{thm:covering-object}(2) that ${\rm Ob}(\mcJ)$ is enveloping. Thus $K$ has an ${\rm Ob}(\mcJ)$-envelope $j:K\to Y$, it follows from Wakamatsu's Lemma (see \cite[Theorem 10.3]{FH} for instance) that there exists a conflation
$$\xymatrix@1@C=45pt{K \ar[r]^{j} & Y \ar[r] & X,}$$
where $X$ is in ${\rm Ob}(\mcI)$. Consider the following pushout diagram
$$\xymatrix@C=30pt@R=30pt{K \ar[r] \ar[d]_{j} & P \ar[r] \ar[d] & A \ar@{=}[d] \\
Y \ar[r]\ar[d] & C \ar[r]^{i} \ar[d] & A \\
X\ar@{=}[r] &X,
}$$
where all rows and columns are conflations. Since $P\in{{\rm Ob}(\mcI)}$, it follows from Lemma \ref{lemma:3.3-1}  that $C\in{{\rm Ob}(\mcI)}$. This shows that $i:C\to A$ is a special ${\rm Ob}(\mcI)$-precover. So $({\rm Ob}(\mcI),{\rm Ob}(\mcJ))$ is a complete cotorsion pair by Salce's Lemma (see \cite[Lemma 2.3]{HZZ} for instance), as desired.
\end{proof}

\section{Enochs Conjecture for object ideals}
Throughout this section, $R$ is an associative ring with identity and $R$-Mod is the category of left $R$-modules.  As stated in Introduction, a hypothetical general positive answer to this question whether each covering ideal in $R$-Mod is closed under direct limits, is sometimes called ``Enochs Conjecture" for objects. In combination this with Proposition \ref{prop:2.1}, a hypothetical general positive answer to this question whether each covering ideal in $R$-Mod is closed under direct limits, is called ``Enochs Conjecture" for morphisms. The aim of this section is to establish the relations between Enochs Conjecture for objects and morphisms. To that end, we need the following lemma.

\begin{lem}\label{lem:direct-summand} The following are ture for any class $\mcX$ in {\rm $R$-Mod}:
\begin{enumerate}

\item If $\mcX$ is covering in {\rm $R$-Mod}, then $\mcX$ is closed under direct summands.

\item If $\mcX$ is enveioping in {\rm $R$-Mod}, then $\mcX$ is closed under direct summands.
\end{enumerate}
\end{lem}
\begin{proof} We only prove (1), and the proof of (2) is dual. Let $M$ be an object in $\mcX$ and $N$ isomorphic to a direct summand of $M$. Then there exist morphisms $i:N\to M$ and $\pi:M\to N$ such that $\pi{i}=1_{N}$. Since $\mcX$ is covering, $N$ has an $\mcX$-cover $f:X\to N$, and therefore there exists $g:M\to X$ such that $\pi=fg$. This yields that $f(gi)=(fg)i={\pi}i=1_{N}$, and hence $f(gi)f=f$. Because $f$ is an $\mcX$-cover, it follows that $(gi)f\cong1_{X}$. So $N\cong X$, as desired.
\end{proof}

Now we state and prove the main result of this section.

\begin{thm} \label{cor1} Let $\mcI$ be an object ideal in {\rm $R$-Mod}. Then $\mcI$ satisfies Enochs Conjecture if and only if ${\rm Ob}(\mcI)$ satisfies Enochs Conjecture.
\end{thm}
\begin{proof} For the ``if" statement, we assume that $\mcI$ is a covering ideal in $R$-Mod. Let $\{f_{i}:M_{i}\to N_{i}\}_{i\in{\mathbb{N}}}$ be a morphism between directed systems of morphisms $\{g_{ij}:M_{i}\to M_{j}\}_{i\leqslant j}$ and $\{h_{ij}:N_{i}\to M_{j}\}_{i\leqslant j}$,
satisfying that $f_{i}\in{\mcI}$ for every $i\in{\mathbb{N}}$. By assumption, each $N_{i}$ has an $\mcI$-cover $\alpha_{i}:X_{i}\to N_{i}$, which implies that there exists $\beta_{i}:M_{i}\to X_{i}$ such that $f_{i}=\alpha_{i}\beta_{i}$ for each $i\in{\mathbb{N}}$. It follows from Proposition \ref{thm:covering-object}(1) that ${\rm Ob}(\mcI)$ is a covering class in $R$-Mod, and therefore each $N_{i}$ has an ${\rm Ob}(\mcI)$-cover $\gamma_{i}:D_{i}\rightarrow N_{i}$. Since $\alpha_{i}$ belongs to the object ideal $\mcI$, it is straightforward to show that there exists a morphism $\mu_{i}:X_{i}\rightarrow D_{i}$ such that $\alpha_{i}=\gamma_{i}\mu_{i}$. Notice that $\alpha_{i}:X_{i}\rightarrow N_{i}$
is an $\mcI$-cover. Then there exists a morphism $\nu_{i}:D_{i}\rightarrow X_{i}$ such that $\gamma_{i}=\alpha_{i}\nu_{i}$, and therefore $\alpha_{i}(\nu_{i}\mu_{i})=(\alpha_{i}\nu_{i})\mu_{i}=\gamma_{i}\mu_{i}=\alpha_{i}$. Since $\alpha_{i}:X_{i}\to N_{i}$ is an $\mcI$-cover, it follows that $\nu_{i}\mu_{i}$ is an isomorphism, and hence $X_{i}$ is isomorphic to a direct summand of $D_{i}$. So $D_{i}\in{{\rm Ob}(\mcI)}$ by Lemma \ref{lem:direct-summand}(1).  Note that ${\rm Ob}(\mcI)$ is closed under direct limits by hypothesis. It follows that $\lim\limits_{\longrightarrow}  X_{i}$ belongs to ${\rm Ob}(\mcI)$. This yields that $\lim\limits_{\longrightarrow} f_{i}=\lim\limits_{\longrightarrow} \alpha_{i}\lim\limits_{\longrightarrow}  \beta_{i}$ belongs to $\mcI$, as desired.

For the ``only if" statement, we assume that ${\rm Ob}(\mcI)$ is a covering class in $R$-Mod. It follows from Proposition \ref{thm:covering-object}(1) that the ideal $\mcI$ is covering. Let $\{X_{i}\}_{i\in{\mathbb{N}}}$ be a family of objects in $\mcX$. Then $1_{X_{i}}\in{\mcI(\mcX)}$ for any $i\in{\mathbb{N}}$. Note that $\mcI$ satisfies Enochs Conjecture by hypothesis. It follows that $\lim\limits_{\longrightarrow}1_{X_{i}}$ belongs to $\mcI$, and therefore $1_{\lim\limits_{\longrightarrow}X_{i}}:\lim\limits_{\longrightarrow}X_{i}\to \lim\limits_{\longrightarrow}X_{i}$ factors through an object $W$ in ${\rm Ob}(\mcI)$. This yields that $\lim\limits_{\longrightarrow}X_{i}$ is isomorphic to a direct summand of $W$. So $\lim\limits_{\longrightarrow}X_{i}$ belongs to ${\rm Ob}(\mcI)$ by Lemma \ref{lem:direct-summand}(1), as desired.
\end{proof}

\vspace{2mm}
As a direct consequence of Theorem \ref{cor1}, we have the following result.
\begin{cor} The class of projective morphisms in $\textrm{R}$-{\rm Mod} satisfies the Enochs Conjecture.
\end{cor}

As another consequence of Theorem \ref{cor1}, we have the following result.
\begin{cor}\label{cor:4.4} Let $\mcI$ and $\mcJ$ be object ideals in $\textrm{R}$-{\rm Mod}, and let $(\mcI, \mcJ)$ be an ideal cotorsion pair. If $\mcJ$ is closed under direct limits, then $\mcI$ satisfies the Enochs Conjecture.
\end{cor}
\begin{proof} The result follows from Theorem \ref{cor1} and \cite[Theorem 3.6]{AST-2018} by noting that ${\rm Ob}(\mcJ)$ is closed under direct limits provided that $\mcJ$ is closed under direct limits.
\end{proof}

%
We end this paper with the following results, which characterize when the validity of Enchos Conjecture for objects of a general ideal $\mcI$ (not necessarily object ideal) implies the validity of Enchos Conjecture of $\mcI$. To that end, we recall from \cite{FH} that an ideal $\mathcal{I}$ in $R$-Mod is \emph{idempotent} if ${\mcI}^{2}=\mcI$. 

\begin{prop}\label{prop4.5} Let $\mcI$ be an idempotent ideal {\rm $R$-Mod}. If ${\rm Ob}(\mcI)$ satisfies Enochs Conjecture in {\rm $R$-Mod}, then so is the ideal $\mcI$.
\end{prop} 
\begin{proof} Assume that $\mcI$ is a covering ideal in {\rm $R$-Mod}. Since $\mcI$ is an idempotent ideal, it follows from \cite[Proposition 10.5]{FH} that $\mcI$ is an object ideal. So the result follows from Theorem \ref{cor1}.
\end{proof}

\begin{prop}\label{prop-last} Let $(\mcI, \mcJ)$ be an ideal cotorsion pair in {\rm $R$-Mod}, and suppose that $\mcJ$ is enveloping and closed under {\rm ME}-extensions. If ${\rm Ob}(\mcI)$ satisfies Enochs Conjecture in {\rm $R$-Mod}, then so is the ideal $\mcI$.
\end{prop}
\begin{proof} We first claim that every $\mcI$-cover of any object $A$ in $R$-Mod is also an ${\rm Ob}(\mcI)$-cover of $A$. Now let $f:D\to A$ be an $\mcI$-cover of $A$. Then $A$ has a conflation
$$\xymatrix@1@C=45pt{K \ar[r] & P \ar[r] & A,}$$
where $P$ is a projective module in $R$-Mod. Since $\mcJ$ is enveloping and closed under {\rm ME}-extensions, it follows from the dual of Lemma 10.1 in \cite{FH} that there exists a conflation $$\xymatrix@1@C=45pt{K \ar[r]^{j} & Y \ar[r] & X,}$$
where $X$ is in ${\rm Ob}(\mcI)$. Consider the following pushout diagram
$$\xymatrix@C=30pt@R=30pt{K \ar[r] \ar[d]_{j} & P \ar[r] \ar[d] & A \ar@{=}[d] \\
Y \ar[r]\ar[d] & C \ar[r]^{i} \ar[d] & A \\
X\ar@{=}[r] &X,
}$$
where all the rows and columns are conflations. Since $P$ and  $X$ belong to ${\rm Ob}(\mcI)$, so is $C$ by Lemma \ref{lemma:3.3-1}. It follows from Lemma \ref{prop:precovering-object}(2) that $i:C\to A$ is also an $\mcI$-precover. Notice that $f:D\to A$ is an $\mcI$-cover of $A$. Then $D$ is isomorphic to a direct summand of $C$, and therefore $D\in{{\rm Ob}(\mcI)}$ by Lemma \ref{lem:direct-summand}(1), as desired.

Next, we assume that $\mcI$ is a covering ideal in $R$-Mod. Let $\{f_{i}:M_{i}\to N_{i}\}_{i\in{\mathbb{N}}}$ be a morphism between directed systems of morphisms $\{g_{ij}:M_{i}\to M_{j}\}_{i\leqslant j}$ and $\{h_{ij}:N_{i}\to M_{j}\}_{i\leqslant j}$,
satisfying that $f_{i}\in{\mcI}$ for every $i\in{\mathbb{N}}$. By assumption, each $N_{i}$ has an $\mcI$-cover $\alpha_{i}:X_{i}\to N_{i}$, which implies that there exists $\beta_{i}:M_{i}\to X_{i}$ such that $f_{i}=\alpha_{i}\beta_{i}$ for each $i\in{\mathbb{N}}$. By the proof above, it follows that ${\rm Ob}(\mcI)$ is a covering class in $\mcA$, and therefore ${\rm Ob}(\mcI)$ is closed under direct limits by hypothesis. Thus $\lim\limits_{\longrightarrow}  X_{i}$ belongs to ${\rm Ob}(\mcI)$. This yields that $\lim\limits_{\longrightarrow} f_{i}=\lim\limits_{\longrightarrow} \alpha_{i}\lim\limits_{\longrightarrow}  \beta_{i}$ belongs to $\mcI$, as desired.
\end{proof}

\bigskip

\renewcommand\refname

\vspace{4mm}
\small

\noindent\textbf{Dandan Sun}\\
School of Mathematical Sciences, Zhejiang University of Technology, Hangzhou 310023, China\\
E-mail: 13515251658@163.com\\[1mm]
\noindent\textbf{Qikai Wang}\\
School of Mathematical Sciences, Zhejiang University of Technology, Hangzhou 310023, China\\
E-mail: qkwang@zjut.edu.cn\\[1mm]
\textbf{Haiyan Zhu}\\
School of Mathematical Sciences, Zhejiang University of Technology, Hangzhou 310023, China\\
E-mails: hyzhu@zjut.edu.cn\\[1mm]
\end{document}